\documentclass[12pt]{amsart}
\usepackage[utf8]{inputenc}

\usepackage{a4}
\usepackage{mathtools}
\usepackage[dvipsnames]{xcolor}
\usepackage{etex}
\usepackage[usenames,dvipsnames]{pstricks} 
\usepackage{epsfig}
\usepackage{graphicx,color}
\usepackage{geometry}
\usepackage[all]{xy}
\usepackage{amssymb,amscd}
\usepackage{cite}
\usepackage{fullpage}
\usepackage{marvosym}
\xyoption{poly}
\usepackage{url}
\usepackage{comment}

\usepackage{tikz}
\usepackage{tikz-cd}
\usetikzlibrary{decorations.pathmorphing}

\theoremstyle{definition}

\newtheorem*{question*}{Question}

\newtheorem*{questions*}{Questions}

\newtheorem*{steps*}{Answer/steps}

\newtheorem*{progress*}{Progress}

\newtheorem*{classification*}{Classification}

\newtheorem*{construction*}{Classification}
\newtheorem*{example*}{Example}

\newtheorem*{remark*}{Remark}
\newtheorem*{remarks*}{Remarks}
\newtheorem*{definition*}{Definition}

\begin{document}

\title{Corrigendum to\\ ``Isomorphism classes of Drinfeld modules over finite fields''}
\author{Valentijn Karemaker, Jeffrey Katen, and Mihran Papikian}
\maketitle

In this note we provide corrections to \cite[Theorem 5.4]{KKP}. The main theorems of \cite{KKP} --Theorem A and B in the introduction-- are valid as stated; in the proof of Theorem B the argument needs to be modified by replacing the erroneous Theorem 5.4 by the theorem of this note. 

References below are to other results in \cite{KKP}. Notation is also as in \emph{loc.~cit.}, of which we recall the following for completeness:\\

\noindent {\bfseries Notation.} Let $\mathbb{F}_q$ be a finite field with $q$ elements. Let $F$ be a function field of transcendence degree $1$ 
over $\mathbb{F}_q$; we assume that $\mathbb{F}_q$ is algebraically closed in $F$. 
Fix a place $\infty$ of $F$. 
Let $A$ be the ring of elements of $F$ regular outside of $\infty$. 
Let $\mathfrak{p}\lhd A$ be a prime, and denote $\mathbb{F}_{\mathfrak{p}} = A/\mathfrak{p}$. 
Let $d=[\mathbb{F}_{\mathfrak{p}}:\mathbb{F}_q]$. 
Let $k\simeq \mathbb{F}_{q^n}$ be a finite extension of $\mathbb{F}_{\mathfrak{p}}$, so $m:= [k:\mathbb{F}_{\mathfrak{p}}]=n/d$ is an integer. 
We consider $k$ as an $A$-field via $\gamma\colon A\to A/\mathfrak{p}\hookrightarrow k$.  

Let $\phi\colon A\to k\{\tau\}$ be a Drinfeld module of rank $r$ over $k$, where $\tau \alpha=\alpha^q\tau$ for all $\alpha\in k$.  
Let $\mathcal{E} := \mathrm{End}_k(\phi)$, $D := \mathcal{E}\otimes_A F$, $\pi:= \tau^n\in \mathcal{E}$, and $\widetilde{F} := F(\pi) \subseteq D$. 
Then $D$ is a central division algebra over $\widetilde{F}$ of dimension $(r/[\widetilde{F}:F])^2$. We restrict to the case where $D$ (or equivalently $\mathcal{E}$) is commutative, which is equivalent to requiring that $r = [\widetilde{F}:F]$, so to requiring that $D = \widetilde{F}$. In this case, $\mathcal{E}$ is an order of $\widetilde{F}$ containing $A[\pi]$. In this erratum, we will restrict to the case where $\mathcal{E} = A[\pi]$.\\

Part (1) of \cite[Theorem 5.4]{KKP} claims the existence of an \emph{action} (given by the map $I \mapsto I * \phi$) of the monoid of fractional ideals of $\mathcal{E}$ up to linear equivalence, on the set of isomorphism classes of Drinfeld modules in the
isogeny class of $\phi$ whose endomorphism ring is the order of an $\mathcal{E}$-ideal. 

This is false: since the endomorphism ring of $I*\phi$ may not coincide with that of $\phi$ (cf. Lemma 4.2), there is no action, unless one restricts to invertible ideals; in this case, the action was decribed previously by Hayes in~\cite{HayesCFT}. 

Instead, we will prove below that the map $I \mapsto I * \phi$ provides a \emph{bijection} between the linear equivalence classes of fractional ideals of $\mathcal{E}$ and the isomorphism classes of Drinfeld modules in the isogeny class of $\phi$.\\

Part (3) of \cite[Theorem 5.4]{KKP} moreover claims that when $\mathcal{E}$ is Gorenstein, the association $I \mapsto I * \phi$ is transitive on the set of all Drinfeld modules whose endomorphism ring is the order of an $\mathcal{E}$-ideal. 

This is false, as is illustrated by the following example, observed by Bergstr{\"o}m and de Vries in \cite[Proposition A.2]{BdV}:\vspace{0.1in}

\noindent {\bfseries Example.}
Let $\mathfrak{p} = T$ and let $\phi$ be the supersingular Drinfeld module of rank $2$ over $k = \mathbb{F}_{q^2}$ defined by $\phi_T = \alpha \tau^2$, where $\alpha \in k^{\times}$. 
 We see that $\widetilde{F}=F(\pi)$ has degree $2$ over~$F$ if and only if $\alpha\not \in \mathbb{F}_q^\times$. Assume this is the case, so that $D = \widetilde{F}$ is commutative and $\mathcal{E}=k[T]$ is the ring of integers (i.e., the maximal order) in $\widetilde{F}$. Now \cite[Theorem 5.4.(3)]{KKP} would imply there is a unique Drinfeld module in the $k$-isogeny class of $\phi$, since all fractional $\mathcal{E}$-ideals are linearly equivalent. However, there are in fact two non-isomorphic Drinfeld modules in this isogeny class, $\phi$ and $\psi$, defined by $\phi_T$ and $\psi_T=\alpha^q \tau^2$ respectively, with the isogeny $\phi\to \psi$ given by $\tau$.\\

In the proof below, we will explain this phenomenon, by showing precisely how non-equivalent ideals give rise to non-isomorphic Drinfeld modules, and hence we will show (in a remark) that in this example there are indeed exactly two isomorphism classes of Drinfeld modules in the $k$-isogeny class of $\phi$.\\
 
The following theorem thus replaces the incorrect \cite[Theorem 5.4]{KKP}. It still implies Corollary 5.5 -- and thus Theorem~B -- of \emph{loc. cit}, while being strictly stronger than Theorem B.
The proof crucially uses ideas from \cite{BKM, LaumonCDV, Waterhouse}, that were not used in the proof of \cite[Theorem 5.4]{KKP}.\\

\noindent {\bfseries Theorem.} Assume that either $k = \mathbb{F}_{\mathfrak{p}}$ or the isogeny class that we consider is ordinary, so that there is a Drinfeld module $\phi$ with $\mathrm{End}_k(\phi) = A[\pi]$. Then the map $I \mapsto I * \phi$ from
the linear equivalences classes of ideals of $A[\pi]$ to the isomorphism classes of Drinfeld modules isogenous to $\phi$ is a bijection.

\begin{proof}
By Lemma 4.1.(1), we may consider the fractional ideals of $\mathcal{E} = A[\pi]$ up to linear equivalence. By Lemma 4.1.(2), the association $I \mapsto I * \phi$ is injective when restricted to kernel ideals. Since $A[\pi]$ is Gorenstein, every $A[\pi]$-ideal is a kernel ideal by Proposition 4.5. Hence, the map $I \mapsto I * \phi$ from linear equivalence classes of $A[\pi]$-ideals is injective.\\ 

It remains to prove surjectivity, i.e., we need to show that every isogeny $u\colon \phi \to \psi$ defined over $k$ arises from some $A[\pi]$-ideal $I$ (up to linear equivalence) in the sense that $\mathbb{H}(\psi) = I \mathbb{H}(\phi)$, cf.~Lemma~3.6. Since $\mathbb{H}(\phi) = \prod_{\mathfrak{l} \lhd A} H_{\mathfrak{l}}(\phi)$ (and similarly for $\psi$), we work locally.\\ 

Let $T_\mathfrak{l}(\phi)$ be the Tate module of $\phi$, where $\mathfrak{l}\lhd A$ is a prime different from $\mathfrak{p}$, and let $H_{\mathfrak{l}}(\phi)=\mathrm{Hom}(T_\mathfrak{l}(\phi), A_\mathfrak{l})$ be its $A_\mathfrak{l}$-dual. 
Since $A[\pi]$ is Gorenstein, by \cite[Theorem~4.9]{GP2}, $T_\mathfrak{l}(\phi)$ is free over $\mathcal{E}_\mathfrak{l}$ of rank $1$, and so is $H_{\mathfrak{l}}(\phi)$ (by a property of Gorenstein rings, cf.~\cite[Definition 4.8]{GP2}. Hence, every $\mathcal{E}_{\mathfrak{l}}$-submodule of $H_{\mathfrak{l}}(\phi)$ is represented by an $\mathcal{E}_{\mathfrak{l}}$-ideal.

Next, let $u\colon \phi\to \psi$ be an isogeny defined over $k$. It induces a natural inclusion 
\begin{equation}\label{eq:incH}
H_\mathfrak{l}(\psi) \subseteq H_\mathfrak{l}(\phi)
\end{equation}
as sublattices. Moreover, by the Tate isomorphism
\begin{equation}\label{eq:Tate}
    \mathrm{Hom}_k(\phi,\psi) \otimes A_{\mathfrak{l}} \xrightarrow{\simeq} \mathrm{Hom}_{A_{\mathfrak{l}}[\pi]}(H_{\mathfrak{l}}(\psi), H_{\mathfrak{l}}(\phi)),
\end{equation}
where the set on the right denotes $(\mathcal{E}_{\mathfrak{l}} =) A_{\mathfrak{l}}[\pi]$-invariant maps, the inclusion of \eqref{eq:incH} realises $H_{\mathfrak{l}}(\psi)$ as an $\mathcal{E}_{\mathfrak{l}}$-submodule of $H_{\mathfrak{l}}(\phi)$.
This implies that
\begin{equation}\label{eq:locallideal}
 H_{\mathfrak{l}}(\psi) = I_{\mathfrak{l}}H_{\mathfrak{l}}(\phi) \text{ for some ideal $I_{\mathfrak{l}} \lhd \mathcal{E}_{\mathfrak{l}}$.}   
\end{equation}
Hence, the above gives us such ideals $I_{\mathfrak{l}}$ for any prime $\mathfrak{l} \neq \mathfrak{p}$ of $A$. Moreover, we have $I_{\mathfrak{l}} = \mathcal{E}_{\mathfrak{l}}$ for all but finitely many $\mathfrak{l}$, since the isogeny $u$ has finite kernel $G$.\\

Now consider $\mathfrak{p} \lhd A$ and the Dieudonn{\'e} module $(H_{\mathfrak{p}}(\phi),f_{\phi,\mathfrak{p}})$; we will drop the $f_{\phi,\mathfrak{p}}$ from the notation when no confusion can arise.  
Let~$F_k$ denote the unramified extension of $F_{\mathfrak{p}}$ with residue field $k$, and $\mathcal{O}_k \simeq k[[\mathfrak{p}]]$ its ring of integers, which is a free $A_{\mathfrak{p}}$-module of rank $m = n/d$. 
Recall that by \cite[\S 2.5]{LaumonCDV}, $H_{\mathfrak{p}}(\phi)$ is a free (left) $\mathcal{O}_k$-module of rank
\begin{equation}\label{eq:OkrankHp}
    \mathrm{rank}_{\mathcal{O}_k}(H_{\mathfrak{p}}(\phi)) = r.
\end{equation}
Morphisms between Dieudonn{\'e} modules are a priori $\mathcal{O}_k$-module morphisms that commute with $f_{\phi,\mathfrak{p}}$, which we will call $\mathcal{O}_k[f_{\phi,\mathfrak{p}}]$-module homomorphisms.

Let $u: \phi \to \psi$ be a $k$-isogeny as before. By the Tate isomorphism at $\mathfrak{p}$, (see \cite[Theorem~2.5.6]{LaumonCDV}) we have
\begin{equation}\label{eq:Tateatp}
 \mathrm{Hom}_k(\phi,\psi) \otimes A_{\mathfrak{p}} \xrightarrow{\simeq} \mathrm{Hom}_{A_{\mathfrak{p}}[\pi]}(H_{\mathfrak{p}}(\psi), H_{\mathfrak{p}}(\phi)).
\end{equation}
Since $A_{\mathfrak{p}}[\pi] = \mathcal{E}_{\mathfrak{p}}$, Equation~\eqref{eq:Tateatp} shows that the sublattice $H_{\mathfrak{p}}(\psi)$ is also a (right) $\mathcal{E}_{\mathfrak{p}}$-submodule of $H_{\mathfrak{p}}(\phi)$; thus, $H_{\mathfrak{p}}(\psi)$ is an $\mathcal{O}_k[f_{\phi,\mathfrak{p}}]\otimes_{A_{\mathfrak{p}}}\mathcal{E}_{\mathfrak{p}}$-module. 

In particular, taking $\psi = \phi$, we get that $H_{\mathfrak{p}}(\phi)$ is itself a (right) $\mathcal{E}_{\mathfrak{p}}$-module, which is free since $A[\pi]$ is Gorenstein, again arguing as in \cite[Theorem 4.9]{GP2}. We now compute its $\mathcal{E}_{\mathfrak{p}}$-rank.
Since $\mathcal{O}_k = kA_{\mathfrak{p}}$ is a free $A_\mathfrak{p}$-module of rank $m$, it follows from Equation~\eqref{eq:OkrankHp} that $\mathrm{rank}_{A_{\mathfrak{p}}}(H_{\mathfrak{p}}(\phi)) = mr$. Further, since $\mathrm{rank}_{A_{\mathfrak{p}}}(\mathcal{E}_{\mathfrak{p}}) = r$, we obtain
\begin{equation}\label{eq:EprankHp}
      \mathrm{rank}_{\mathcal{E}_{\mathfrak{p}}}(H_{\mathfrak{p}}(\phi)) = m.
\end{equation}
When $m=1$, or equivalently when $k=\mathbb{F}_\mathfrak{p}$, we get that $\mathcal{O}_k \simeq A_{\mathfrak{p}}$, hence $\mathcal{O}_k[f_{\phi,\mathfrak{p}}] \otimes_{A_{\mathfrak{p}}} \mathcal{E}_{\mathfrak{p}} \simeq \mathcal{E}_{\mathfrak{p}}$, and every $\mathcal{E}_{\mathfrak{p}}$-submodule of $H_{\mathfrak{p}}(\phi)$ is given by an ideal. So in this case, for an isogeny $u\colon \phi \to \psi$ we get 
\begin{equation}\label{eq:localpideal}
    H_{\mathfrak{p}}(\psi) = I_{\mathfrak{p}} H_{\mathfrak{p}}(\phi).
\end{equation} 
as for $\mathfrak{l}$ above. 
We then note that the local ideals obtained in \eqref{eq:localpideal} and \eqref{eq:locallideal} for all $\mathfrak{l} \neq \mathfrak{p}$ are the localisations of a global lattice, which is closed under the action of $\mathcal{E}$, and therefore is a global ideal $I$. This ideal satisfies $\mathbb{H}(\psi) = I \mathbb{H}(\phi)$ by construction, as required.\\

Suppose now that $m > 1$. We do not immediately impose that $\phi$ is ordinary, but below we will indicate precisely when this assumption is needed to conclude the proof.\\
Recall that a priori the endomorphism algebra $D = \widetilde{F} = F(\pi)$ of $\phi$ satisfies
\begin{equation}\label{eq:Dsplit}
D_{\mathfrak{p}} := D \otimes_F F_{\mathfrak{p}} = \bigoplus_{\nu \vert \mathfrak{p}} D_{\nu},
\end{equation}
where $\nu$ runs over all places of the maximal $A$-order $B$ in $D$ which lie over $\mathfrak{p}$. 
This splitting thus descends to~$B$:
\begin{equation}\label{eq:Bsplit}
    B_{\mathfrak{p}} = \bigoplus_{\nu \vert \mathfrak{p}} B_{\nu}.
\end{equation}
Here each $B_{\nu}$ is a local ring with unique maximal ideal $\mathfrak{p}_{B_{\nu}}$. Among $\nu \vert \mathfrak{p}$, we distinguish the place $\widetilde{\mathfrak{p}}$, being the unique place of $D = \widetilde{F}$ lying over the place $(\pi)$ of $K = \mathbb{F}_q(\pi)$.

The endomorphism ring $\mathcal{E} = \mathrm{End}_k(\phi) = A[\pi]$ is an order contained in~$B$, but not necessarily equal to it. Hence, the splitting of~\eqref{eq:Bsplit} may not descend to $\mathcal{E}$. But since $\mathcal{E}_{\mathfrak{p}}$ is complete and semilocal, we do obtain the splitting
\begin{equation}\label{eq:Esplit}
    \mathcal{E}_{\mathfrak{p}} = \bigoplus_{\mathfrak{P} \vert \mathfrak{p}} \mathcal{E}_{\mathfrak{P}},
\end{equation}
where the direct sum runs over the maximal ideals $\mathfrak{P} = \mathfrak{p}_{B_{\nu}} \cap \mathcal{E}_{\mathfrak{p}}$ of $\mathcal{E}_{\mathfrak{p}}$; all maximal ideals of $\mathcal{E}_{\mathfrak{p}}$ arise as such intersections, but those for different $\nu$ may coincide (cf.~\cite[Proof of Theorem~7.4]{Waterhouse} or \cite[\S 4.3.1]{BKM}). This splitting in turn induces an alternative splitting $D_{\mathfrak{p}} = \oplus_{\mathfrak{P} \vert \mathfrak{p}} D_{\mathfrak{P}}$.\\
 
Let
$\widetilde{\mathfrak{P}} = \widetilde{\mathfrak{p}} \cap \mathcal{E}_{\mathfrak{p}}$. By Corollary~2.9 (or alternatively by \cite[Corollary, p.~164]{JKYu}), $\mathcal{E} = A[\pi]$ is locally maximal at $\pi$, i.e., $\mathcal{E}_{\widetilde{\mathfrak{P}}} = B_{\widetilde{\mathfrak{p}}}$ is maximal and $D_{\widetilde{\mathfrak{P}}} = D_{\widetilde{\mathfrak{p}}}$.
We obtain
\begin{equation}
\begin{split}
    \mathcal{E}_{\mathfrak{p}} = \mathcal{E}_{\widetilde{\mathfrak{P}}} \oplus \left( \bigoplus_{\mathfrak{P} \neq \widetilde{\mathfrak{P}}} \mathcal{E}_{\mathfrak{P}} \right) & = B_{\widetilde{\mathfrak{p}}} \oplus \mathcal{E}'_{\mathfrak{p}}, \\
   D_{\mathfrak{p}} = D_{\widetilde{\mathfrak{P}}} \oplus \left( \bigoplus_{\mathfrak{P} \neq \widetilde{\mathfrak{P}}} D_{\mathfrak{P}} \right) & = D_{\widetilde{\mathfrak{p}}} \oplus D'_{\mathfrak{p}}.
\end{split}
\end{equation}

For the Dieudonn{\'e} module $H_{\mathfrak{p}}(\phi)$, as in Proposition 4.5, cf.~\cite[Proposition~2.4.6]{LaumonCDV}, we have a decomposition
\begin{equation}
    H_{\mathfrak{p}}(\phi) = H_{\mathfrak{p}}^c(\phi) \oplus H_{\mathfrak{p}}^{\mathrm{\acute{e}t}}(\phi)
\end{equation}
into its connected component $H_{\mathfrak{p}}^c(\phi)$ (of positive slope) and its maximal {\'e}tale quotient $H_{\mathfrak{p}}^{\mathrm{\acute{e}t}}(\phi)$ (of slope zero). There is a corresponding decomposition of the Frobenius
\begin{equation}
    f_{\phi,\mathfrak{p}} = f_{\phi,\widetilde{\mathfrak{p}}} \oplus f'_{\phi,\mathfrak{p}} = \bigoplus_{\mathfrak{P} \vert \mathfrak{p}}f_{\phi,\mathfrak{P}}
\end{equation}
such that
\begin{equation}
     \mathcal{E}_{\widetilde{\mathfrak{P}}} = B_{\widetilde{\mathfrak{p}}} \simeq \mathrm{End}(H_{\mathfrak{p}}^c(\phi), f_{\phi, \widetilde{\mathfrak{p}}}), \qquad \mathcal{E}'_{\mathfrak{p}} \simeq \mathrm{End}(H_{\mathfrak{p}}^{\mathrm{\acute{e}t}}(\phi), f'_{\phi,\mathfrak{p}}).
\end{equation}
Moreover, $H_{\mathfrak{p}}^c(\phi)$ is a free $B_{\widetilde{\mathfrak{p}}}$-module and $H_{\mathfrak{p}}^{\mathrm{\acute{e}t}}(\phi)$ is a free $\mathcal{E}'_{\mathfrak{p}}$-module.\\

Using the decompositions above, our task is now to describe $(\mathcal{O}_k[f_{\phi,\widetilde{\mathfrak{p}}}] \otimes_{A_{\mathfrak{p}}} B_{\widetilde{\mathfrak{p}}})$-submodules of $H^c_{\mathfrak{p}}(\phi)$, resp.~$(\mathcal{O}_k[f'_{\phi,\mathfrak{p}}] \otimes_{A_{\mathfrak{p}}} \mathcal{E}'_{\mathfrak{p}})$-submodules of $H^{\mathrm{\acute{e}t}}_{\mathfrak{p}}(\phi)$, as ideals, viz.~$I_{\widetilde{\mathfrak{p}}}H^c_{\mathfrak{p}}(\phi)$, resp.~$I'_{\mathfrak{p}} H^{\mathrm{\acute{e}t}}_{\mathfrak{p}}(\phi)$.
Using again that $\mathcal{O}_k = k A_{\mathfrak{p}}$ is a free $A_{\mathfrak{p}}$-module of rank $m$, we see that 
\begin{equation}
\mathcal{O}_k \otimes_{A_{\mathfrak{p}}} \mathcal{E}_{\mathfrak{P}} = \bigoplus_{i=1}^{g_{\mathfrak{P}}} k \mathcal{E}_{\mathfrak{P}}
\end{equation}
for any $\mathfrak{P} \vert \mathfrak{p}$, where $g_{\mathfrak{P}} = \mathrm{gcd}(f(\mathfrak{P}/\mathfrak{p}),m)$ and where $f(\mathfrak{P}/\mathfrak{p})$ is the residue degree of $\mathfrak{P}$. 
Thus, for each $\mathfrak{P} \vert \mathfrak{p}$, the $(\mathcal{O}_k[f_{\phi,\mathfrak{P}}] \otimes_{A_{\mathfrak{p}}} \mathcal{E}_{\mathfrak{P}})$-submodules are $f_{\phi,\mathfrak{P}}$-stable modules of the form $I_1 \oplus \ldots \oplus I_{g_{\mathfrak{P}}}$ where each $I_i$ is a fractional ideal in its respective summand $k \mathcal{E}_{\mathfrak{P}}$.\\

To check stability under $f_{\phi,\mathfrak{p}}$, we must first represent its action on $H_{\mathfrak{p}}(\phi)$, which we do rationally and locally for each place $\mathfrak{P} \vert \mathfrak{p}$, closely following \cite[Theorems~5.1 and~7.4]{Waterhouse} and \cite[Theorem~2.5.6]{LaumonCDV}.
To represent $V_\mathfrak{p}(\phi):= H_\mathfrak{p}(\phi)\otimes_{A_\mathfrak{p}}F_\mathfrak{p}$ itself, we consider the decomposition 
\begin{equation}
    F_k \otimes_{F_{\mathfrak{p}}} D_{\mathfrak{p}} \simeq \bigoplus_{\mathfrak{P} \vert \mathfrak{p}} (F_k \otimes_{F_{\mathfrak{p}}} D_{\mathfrak{P}}).
\end{equation}
By construction, this algebra has a left $F_k$-action and a right $D_{\mathfrak{p}}$-action; recall that $F_k = \mathrm{Frac}(\mathcal{O}_k)$ is the unramified extension of $F_{\mathfrak{p}}$ of degree $m$ with residue field $k$.\\

First we consider $\widetilde{\mathfrak{P}} \vert \mathfrak{p}$. Here the local endomorphism ring $\mathcal{E}_{\widetilde{\mathfrak{P}}} = B_{\widetilde{\mathfrak{p}}}$ is maximal by the above, so we may argue as in \cite[Theorem~5.1]{Waterhouse}:
The field $F_k\cap D_{\widetilde{\mathfrak{P}}}$ is an unramified extension of $F_{\mathfrak{p}}$ of degree $g_{\widetilde{\mathfrak{P}}}$, and the field $F_k D_{\widetilde{\mathfrak{P}}}$ is an unramified extension of $D_{\widetilde{\mathfrak{P}}}$ of degree $\mathrm{lcm}(m, f(\widetilde{\mathfrak{P}}/\mathfrak{p}))/f(\widetilde{\mathfrak{P}}/\mathfrak{p}) = \frac{m \cdot f(\widetilde{\mathfrak{P}}/\mathfrak{p})}{g_{\widetilde{\mathfrak{P}}}\cdot f(\widetilde{\mathfrak{P}}/\mathfrak{p})} = \frac{m}{g_{\widetilde{\mathfrak{P}}}}$.
Hence, we see that $F_k\otimes_{F_\mathfrak{p}} D_{\widetilde{\mathfrak{P}}}$ is a direct sum of $g_{\widetilde{\mathfrak{P}}}$ copies of the compositum $F_kD_{\widetilde{\mathfrak{P}}}$:
\begin{equation}
    F_k \otimes_{F_{\mathfrak{p}}} D_{\widetilde{\mathfrak{P}}} \xrightarrow{\simeq} F_kD_{\widetilde{\mathfrak{P}}} \oplus \ldots \oplus F_kD_{\widetilde{\mathfrak{P}}}.
\end{equation}
Let $\sigma$ be the Frobenius of $F_k$ over $F_{\mathfrak{p}}$. The map giving the above identification is 
\begin{equation}\label{eq1}
    \omega\otimes \beta\longmapsto (\omega\beta, \sigma(\omega)\beta, \dots, \sigma^{g_{\widetilde{\mathfrak{P}}}-1}(\omega)\beta).     
\end{equation}
An element $\lambda \in F_k$ acts on the direct sum by the diagonal matrix 
\begin{equation}
    \mathrm{diag}(\lambda, \sigma(\lambda), \ldots, \sigma^{g_{\widetilde{\mathfrak{P}}}-1}(\lambda)).
\end{equation}
Furthermore, $\sigma$ acting on the $F_k$-factor of the tensor product takes $\omega\otimes \beta$ 
to 
\begin{equation}\label{eq2}
(\sigma(\omega)\beta, \sigma^2(\omega)\beta, \dots, \sigma^{g_{\widetilde{\mathfrak{P}}}}(\omega)\beta).
\end{equation}
Note that $\sigma$ has order $m$, hence $\sigma^{g_{\widetilde{\mathfrak{P}}}}$ has order $m/g_{\widetilde{\mathfrak{P}}}$, which is the residue degree of $F_k D_{\widetilde{\mathfrak{P}}}$ over $D_{\widetilde{\mathfrak{P}}}$. 
Also note that the Frobenius $\mathrm{Fr}_{\widetilde{\mathfrak{P}}}$ of $F_kD_{\widetilde{\mathfrak{P}}}$ over $D_{\widetilde{\mathfrak{P}}}$ fixes $D_{\widetilde{\mathfrak{P}}}$. 
Thus, 
$\sigma^{g_{\widetilde{\mathfrak{P}}}}(\omega)\beta = \mathrm{Fr}_{\widetilde{\mathfrak{P}}}(\omega\beta)$, so $\sigma$ acts by a cyclic permutation followed by $\mathrm{Fr}_{\widetilde{\mathfrak{P}}}$ in the last place. 
We now use \cite[(4.1.8)]{DM} evaluated on $\pi = \tau^n$ to obtain the equality
\begin{equation}
    f(\widetilde{\mathfrak{P}}/\mathfrak{p})\cdot \mathrm{ord}_{\widetilde{\mathfrak{P}}}(\pi) = \frac{n\cdot f(\widetilde{\mathfrak{P}}/\mathfrak{p}) e(\widetilde{\mathfrak{P}}/\mathfrak{p})}{d\cdot H(\phi)} = \frac{n\cdot [D_{\widetilde{\mathfrak{P}}}:F_\mathfrak{p}]}{d\cdot H(\phi)} = m,
\end{equation}
where $e(\widetilde{\mathfrak{P}}/\mathfrak{p})$ is the ramification index of $\widetilde{\mathfrak{P}}$ and where the last equality follows from our assumption that $D$ is commutative. 
Hence, $m/g_{\widetilde{\mathfrak{P}}}$ divides $\mathrm{ord}_{\widetilde{\mathfrak{P}}}(\pi)$. 
 By local class field theory, if $K_s$ is the unramified extension of a local field $K$ of degree $s$, then 
 \begin{equation}
 \mathrm{Nr}_{K_s/K}(K_s^\times) = \{a\in K^\times\ \colon\  s\vert\mathrm{ord}_K(a)\}=\mathcal{O}_K^\times\cdot {\varpi_K}^{s\mathbb{Z}}.
 \end{equation}
 Hence, the element $\pi\in D_{\widetilde{\mathfrak{P}}}$ is the norm of some element $\alpha\in F_kD_{\widetilde{\mathfrak{P}}}$: 
 \begin{equation}
 \pi=\mathrm{Nr}_{F_kD_{\widetilde{\mathfrak{P}}}/D_{\widetilde{\mathfrak{P}}}}(\alpha) =\alpha\cdot \mathrm{Fr}_{\widetilde{\mathfrak{P}}}(\alpha)\cdots\mathrm{Fr}_{\widetilde{\mathfrak{P}}}^{\frac{m}{g_{\widetilde{\mathfrak{P}}}}-1}(\alpha).
 \end{equation}
 Let $u=(1, 1, \dots, \alpha)\in \bigoplus F_kD_{\widetilde{\mathfrak{P}}}$ and define $f_{\phi, \widetilde{\mathfrak{P}}}=u\sigma$. 
 Then $f_{\phi, \widetilde{\mathfrak{P}}}\lambda=\lambda^\sigma f_{\phi, \widetilde{\mathfrak{P}}}$ for all $\lambda\in F_k$ and 
 \begin{align*}
 f_{\phi, \widetilde{\mathfrak{P}}}^{m} =(u\sigma)^{m} =u u^{\sigma} \cdots u^{\sigma^{m}-1} \sigma^{m} =u u^{\sigma} \cdots u^{\sigma^{m}-1}. 
 \end{align*}
 Now 
 \begin{align*}
 & u^\sigma=(1,\dots, \alpha, 1), \dots, u^{\sigma^{{g_{\widetilde{\mathfrak{P}}}}-1}}=(\alpha, 1, \dots, 1), u^{\sigma^{{g_{\widetilde{\mathfrak{P}}}}}}=(1, 1, \dots, \mathrm{Fr}_{\widetilde{\mathfrak{P}}}(\alpha)),\\ 
 & u^{\sigma^{{g_{\widetilde{\mathfrak{P}}}}+1}}=(1,\dots,  \mathrm{Fr}_{\widetilde{\mathfrak{P}}}(\alpha), 1), \dots, u^{\sigma^{{2g_{\widetilde{\mathfrak{P}}}}-1}}=(\mathrm{Fr}_{\widetilde{\mathfrak{P}}}(\alpha), 1, \dots, 1), u^{\sigma^{{2g_{\widetilde{\mathfrak{P}}}}}}=(1, 1, \dots, \mathrm{Fr}_{\widetilde{\mathfrak{P}}}^2(\alpha)), \\ 
 & \dots u^{\sigma^{m}-1}  = (\mathrm{Fr}_{\widetilde{\mathfrak{P}}}^{\frac{m}{g_{\widetilde{\mathfrak{P}}}}-1}(\alpha), 1, \dots, 1).
 \end{align*}
 We conclude that 
 \begin{equation}
 f_{\phi, \widetilde{\mathfrak{P}}}^{m} = u u^{\sigma} \cdots u^{\sigma^{m}-1} =\left(\mathrm{Nr}_{F_kD_{\widetilde{\mathfrak{P}}}/D_{\widetilde{\mathfrak{P}}}}(\alpha), \cdots, \mathrm{Nr}_{F_kD_{\widetilde{\mathfrak{P}}}/D_{\widetilde{\mathfrak{P}}}}(\alpha)\right) = (\pi, \dots, \pi) =\pi. 
 \end{equation}
 Thus, we have constructed the algebra acting on $\oplus F_kD_{\widetilde{\mathfrak{P}}}$, providing a representation of $f_{\phi,\widetilde{\mathfrak{P}}}$. As this space has the right dimension, it is 
 isomorphic to the $\widetilde{\mathfrak{P}}$-component of $V_\mathfrak{p}(\phi):= H_\mathfrak{p}(\phi)\otimes_{A_\mathfrak{p}}F_\mathfrak{p}$. On this representation, we may now check $f_{\phi,\widetilde{\mathfrak{P}}}$-stability of our fractional ideals.\\

The above simplifies significantly when we use the assumption that $\phi$ is ordinary.
Indeed, then we have $H(\phi) = 1$. 
This implies that $f(\widetilde{\mathfrak{P}}/\mathfrak{p}) = 1$ since ($D_{\widetilde{\mathfrak{P}}} = D_{\widetilde{\mathfrak{p}}} =) \widetilde{F}_{\widetilde{\mathfrak{p}}}= F_{\mathfrak{p}}$. So $g_{\widetilde{\mathfrak{P}}} = 1$ and $F_k \otimes_{F_{\mathfrak{p}}} D_{\widetilde{\mathfrak{P}}} \simeq F_k D_{\widetilde{\mathfrak{P}}} = F_k$, on which $f_{\phi, \widetilde{\mathfrak{p}}} = \alpha \sigma$ acts. 
From $D_{\widetilde{\mathfrak{P}}} = F_{\mathfrak{p}}$ we also get that its maximal order $\mathcal{E}_{\widetilde{\mathfrak{P}}}$ equals $A_{\mathfrak{p}}$, and so $\mathcal{O}_k\otimes_{A_{\mathfrak{p}}} \mathcal{E}_{\widetilde{\mathfrak{P}}} \simeq \mathcal{O}_k \simeq k[[\mathfrak{p}]]$. The submodules we are considering therefore simplify to $\mathcal{O}_k$-fractional ideals $I_{\widetilde{\mathfrak{P}}}$ which are of the form $\mathcal{O}_k\langle \mathfrak{p}^\epsilon \rangle$ for some $\epsilon \geq 0$, i.e. they are generated by $\mathfrak{p}^{\epsilon}$ as $\mathcal{O}_k$-modules. These ideals are visibly stable under $f_{\phi, \widetilde{\mathfrak{P}}}$. 
Hence, we find that every 
$(\mathcal{O}_k[f_{\phi,\widetilde{\mathfrak{P}}}] \otimes_{A_{\mathfrak{p}}} \mathcal{E}_{\widetilde{\mathfrak{P}}})$-submodule of $H^c_{\mathfrak{p}}(\phi)$ is given by $I_{\widetilde{\mathfrak{P}}}H^c_{\mathfrak{p}}(\phi)$, as required.\\

\noindent {\bfseries Lemma.} When $\phi$ is not ordinary, so $g_{\widetilde{\mathfrak{P}}} > 1$, there are ${g_{\widetilde{\mathfrak{P}}} \mathrm{ord}_{\widetilde{\mathfrak{P}}}(\pi)/m+1 \choose g_{\widetilde{\mathfrak{P}}}-1}$ inequivalent $(\mathcal{O}_k[f_{\phi,\widetilde{\mathfrak{p}}}] \otimes_{A_{\mathfrak{p}}} \mathcal{E}_{\widetilde{\mathfrak{P}}})$-submodules of $H_{\mathfrak{p}}^c(\phi)$.

\noindent {\bfseries Proof of lemma.} The arguments of \cite[Theorem~5.1]{Waterhouse} still apply. These show that since $\mathcal{O}_k \otimes_{A_{\mathfrak{p}}} \mathcal{E}_{\widetilde{\mathfrak{P}}} \simeq \oplus k\mathcal{E}_{\widetilde{\mathfrak{P}}}$ is the maximal order in $F_k \otimes D_{\widetilde{\mathfrak{P}}} \simeq \oplus F_k \widetilde{F}_{\widetilde{\mathfrak{p}}}$, 
 the $(\mathcal{O}_k \otimes_{A_{\mathfrak{p}}} \mathcal{E}_{\widetilde{\mathfrak{P}}})$-submodules $I_1 \oplus \ldots \oplus I_{g_{\widetilde{\mathfrak{P}}}}$ are of the form $(\oplus k\mathcal{E}_{\widetilde{\mathfrak{P}}})\langle \widetilde{\mathfrak{P}}^{\epsilon_1}, \ldots, \widetilde{\mathfrak{P}}^{\epsilon_{g_{\widetilde{\mathfrak{P}}}}} \rangle$. 
Invariance under $f_{\phi, \widetilde{\mathfrak{p}}}$ then translates into the condition 
\[
\epsilon_1 \leq \epsilon_2 \leq \ldots \leq \epsilon_{g_{\widetilde{\mathfrak{P}}}} \leq \epsilon_1 + \mathrm{ord}_{\widetilde{\mathfrak{P}}}(\alpha),
\]
 where $\mathrm{ord}_{\widetilde{\mathfrak{P}}}$ is the extension of the normalised valuation on $\widetilde{F}_{\widetilde{\mathfrak{p}}}$ to the unramified extension $F_k \widetilde{F}_{\widetilde{\mathfrak{p}}}$. 
 We note that 
 \[
 \mathrm{ord}_{\widetilde{\mathfrak{P}}}(\alpha)=\mathrm{ord}_{\widetilde{\mathfrak{P}}}(\pi)/[F_k \widetilde{F}_{\widetilde{\mathfrak{p}}} : \widetilde{F}_{\widetilde{\mathfrak{p}}}]=\frac{\mathrm{ord}_{\widetilde{\mathfrak{P}}}(\pi) \cdot g_{\widetilde{\mathfrak{P}}}}{m}. 
 \]
 Since $\mathcal{E}_{\widetilde{\mathfrak{P}}}$ acts diagonally on $\oplus F_k D_{\widetilde{\mathfrak{P}}}$, after scaling by $\widetilde{\mathfrak{P}}^{-\epsilon_1}$, we may assume $\epsilon_1=0$. 
 Hence, the number of $f_{\phi, \widetilde{\mathfrak{p}}}$-stable modules up to scaling is the number of integers 
 \[ 
 0\leq  n_1\leq  \cdots\leq n_{g-1}\leq  \frac{\mathrm{ord}_{\widetilde{\mathfrak{P}}}(\pi) \cdot g_{\widetilde{\mathfrak{P}}}}{m}.  
 \]
 This is exactly the binomial coefficient ${g_{\widetilde{\mathfrak{P}}} \mathrm{ord}_{\widetilde{\mathfrak{P}}}(\pi)/m+1 \choose g_{\widetilde{\mathfrak{P}}}-1}$. \qed

\noindent {\bfseries Remark.} The lemma also explains the above-mentioned example of \cite[Proposition A.2]{BdV}. Recall that in this case $\phi$ is supersingular of rank $2$ over $k = \mathbb{F}_{q^2}$ defined by $\phi_T = \alpha \tau^2$, where $\alpha \in k^{\times}$ is assumed to satisfy $\alpha\not \in \mathbb{F}_q^\times$ so that $\widetilde{F}=F(\pi)$ has degree $2$ over~$F$, and $\mathfrak{p} = T$. Then $D = \widetilde{F}$ is commutative and $\mathcal{E}=k[T]$ is the ring of integers (i.e., the maximal order~$B$) in $\widetilde{F}$. 
 Since $\widetilde{\mathfrak{p}}$ is the only prime above $\mathfrak{p}$ (cf.~\cite[Corollary~4.1.14.(1)]{DM}), it suffices to determine the $(\mathcal{O}_k[f_{\phi,\widetilde{\mathfrak{p}}}] \otimes_{A_{\mathfrak{p}}} \mathcal{E}_{\widetilde{\mathfrak{p}}})$-modules  in $H_{\mathfrak{p}}(\phi) = H^c_{\mathfrak{p}}(\phi)$.
 We have $f(\widetilde{\mathfrak{p}}/\mathfrak{p})=2$ and $m=[k:\mathbb{F}_T]=2$, so $g_{\widetilde{\mathfrak{p}}}=2$. Note that $T$ remains inert in $\widetilde{F}$, so $\mathrm{ord}_{\widetilde{\mathfrak{p}}}(\pi)=\mathrm{ord}_T(T)=1$. By the above, we have ${2 \choose 1}=2$ inequivalent $(\mathcal{O}_k[f_{\phi,\widetilde{\mathfrak{p}}}] \otimes_{A_{\mathfrak{p}}} \mathcal{E}_{\widetilde{\mathfrak{p}}})$-modules  in $H_{\mathfrak{p}}(\phi)$. In the $k$-isogeny class of $\phi$ we correspondingly have exactly two Drinfeld modules, $\phi$ and $\psi$, defined by $\phi_T$ and $\psi_T=\alpha^q \tau^2$ respectively, with the isogeny $\phi\to \psi$ given by $\tau$. This is in accordance with \cite[Proposition A.2]{BdV}.\\ 

Next, we return to the main proof and consider each place $\mathfrak{P}\vert \mathfrak{p}$ such that $\mathfrak{P}\neq \widetilde{\mathfrak{P}}$. Assuming $\phi$ is ordinary, we have $\mathrm{ord}_\mathfrak{P}(\pi)=0$, so $\pi$ is a unit in each $\mathcal{E}_{\mathfrak{P}}$. 
 By a direct adaptation of \cite[Proposition~7.3]{Waterhouse}, which applies Hensel's lemma to the norm form $\mathrm{Nr}_{(\mathcal{O}_k \otimes \mathcal{E}_{\mathfrak{P})}/\mathcal{E}_{\mathfrak{P}}}$, we find that $\pi \in \mathcal{E}_{\mathfrak{P}}$ is the norm of some element $\beta \in \mathcal{O}_k \otimes_{A_{\mathfrak{p}}} \mathcal{E}_{\mathfrak{P}}$, which is thus itself a unit: 
 \begin{equation}
 \pi=\mathrm{Nr}_{(\mathcal{O}_k \otimes \mathcal{E}_{\mathfrak{P})}/\mathcal{E}_{\mathfrak{P}}}(\beta).
 \end{equation}
 As in \cite[Theorem~7.4]{Waterhouse}, let $f_{\phi, \mathfrak{P}}=\beta \sigma$ act on $F_k \otimes_{F_{\mathfrak{p}}} D_{\mathfrak{P}}$, where as before $\sigma$ denotes the Frobenius of $F_k$ over $F_{\mathfrak{p}}$.
 As for $\widetilde{\mathfrak{P}}$, one checks that $f_{\phi,\mathfrak{P}} \lambda = \lambda^{\sigma} f_{\phi,\mathfrak{P}}$ for all $\lambda \in F_k$, and that $f_{\phi,\mathfrak{P}}^m = \pi$.
 Defining $v = \oplus_{\mathfrak{P} \neq \widetilde{\mathfrak{P}}} \beta$ and $f'_{\phi, \mathfrak{p}} = \oplus_{\mathfrak{P} \neq \widetilde{\mathfrak{P}}} f_{\phi, \mathfrak{P}} = v\sigma$ thus provides a suitable representation of the action of Frobenius on $H_{\mathfrak{p}}^{\mathrm{\acute{e}t}}(\phi) \otimes_{A_{\mathfrak{p}}} F_{\mathfrak{p}}$.\\ 

We now argue as in \cite[Proposition~4.16]{BKM} to describe the $(\mathcal{O}_k[f_{\phi,\mathfrak{P}}] \otimes_{A_{\mathfrak{p}}} \mathcal{E}_{\mathfrak{P}})$-modules for all $\mathfrak{P} \neq \widetilde{\mathfrak{P}}$. Recall that $k \simeq \mathbb{F}_{q^n}$ and $\mathbb{F}_{\mathfrak{p}} \simeq \mathbb{F}_{q^d}$, where $m = n/d$ is an integer. Consider the Galois group $G = \mathrm{Gal}(F_k/F_{\mathfrak{p}}) \simeq \mathrm{Gal}(k/\mathbb{F}_{\mathfrak{p}})$, which is cyclic of order $m$, and has generator $\sigma$ which on $k$ acts as $x \mapsto x^{q^d}$. Then $G$ acts also on $\mathcal{O}_k = kA_{\mathfrak{p}} \simeq k[[\mathfrak{p}]]$, and hence on $\mathcal{O}_k \otimes_{A_{\mathfrak{p}}} \mathcal{E}_{\mathfrak{P}}$, in such a way that $\mathcal{O}_k^G = A_{\mathfrak{p}}$ and hence $(\mathcal{O}_k \otimes_{A_{\mathfrak{p}}} \mathcal{E}_{\mathfrak{P}})^G = \mathcal{E}_{\mathfrak{P}}$.

We claim that $(\mathcal{O}_k \otimes_{A_{\mathfrak{p}}} \mathcal{E}_{\mathfrak{P}})/\mathcal{E}_{\mathfrak{P}}$ is a $G$-Galois extension of rings. By \cite[Theorem 1.6.(ii')]{Greither}, this holds if and only if there exist $n \in \mathbb{N}$ and $x_1,\ldots, x_n, y_1, \ldots, y_n \in \mathcal{O}_k \otimes_{A_{\mathfrak{p}}} \mathcal{E}_{\mathfrak{P}}$ such that 
\[
\sum_{i=1}^n x_i \rho(y_i) = \begin{cases} 1 & \text{if } \rho = 1; \\ 0 & \text{ if } \rho \neq 1.   
\end{cases}
\]
Since $G$ acts trivially on $\mathcal{E}_{\mathfrak{P}}$, it suffices to choose $x_i, y_i \in \mathcal{O}_k$. Every element $1 \neq x \in k$ satisfies 
\[
1+x+\ldots+x^{q^n-2} = \sum_{i=1}^{q^n-1} x^{i-1}=0.
\]
So choose such an $x \neq 1$, let $n = q^n-1$, and choose $x_i = \frac{x^{-i}}{q^n - 1}$ and $y_i = x^i$ for any $i = 1, \ldots, q^n-1$. Then
\[
\sum_{i=1}^{n} x_i y_i = \sum_{i=1}^{q^n-1} \frac{x^{-i}}{q^n - 1} x^i = 1
\]
and, writing $\rho= \sigma^j$ for some $1 \leq j \leq m-1$,
\[
\sum_{i=1}^n x_i \rho(y_i) = \frac{1}{q^n-1} \sum_{i=1}^{q^n-1} (x^{q^{dj}-1})^i = 0.
\]
This proves the claim.

Consider the descent datum $\Phi = \{\Phi_{\rho} = \rho\}_{\rho \in G}$, where we view each $\rho$ as an $\mathcal{E}_{\mathfrak{P}}$-automorphism of $\mathcal{O}_k \otimes_{A_{\mathfrak{p}}} \mathcal{E}_{\mathfrak{P}}$ by the discussion above. Then by~\cite[Theorem~7.1]{Greither}, for any $(\mathcal{O}_k \otimes_{A_{\mathfrak{p}}} \mathcal{E}_{\mathfrak{P}})$-module $J$ we obtain by descent an $\mathcal{E}_{\mathfrak{P}}$-module 
\[
J^{\Phi} = \{ j \in J: \Phi_{\rho}(j) = \rho(j) = j \text{ for all } \rho \in G \} = J^G,
\]
such that there is a module isomorphism $J \simeq (\mathcal{O}_k \otimes_{A_{\mathfrak{p}}} \mathcal{E}_{\mathfrak{P}}) \otimes_{\mathcal{E}_{\mathfrak{P}}} J^{\Phi}$. 

Now consider a $(\mathcal{O}_k [f_{\phi,\mathfrak{P}}]\otimes_{A_{\mathfrak{p}}} \mathcal{E}_{\mathfrak{P}})$-module $J \subseteq F_k \otimes D_{\mathfrak{P}}$.  
By the above, $f_{\phi, \mathfrak{P}} = \beta \sigma$ where $\beta$ is a unit. Hence, $\beta J = J$ and since $G = \langle \sigma \rangle$, stability under $f_{\phi, \mathfrak{P}}$ is equivalent to stability under $G$. 
It follows that $J^G = J \cap D_{\mathfrak{P}}$, where by slight abuse of notation we identify $D_{\mathfrak{P}}$ with its (diagonal) image in $F_k \otimes D_{\mathfrak{P}}$. We get that 
\[
J \simeq (\mathcal{O}_k \otimes_{A_{\mathfrak{p}}} \mathcal{E}_{\mathfrak{P}}) \otimes_{\mathcal{E}_{\mathfrak{P}}} J^{G} \simeq \mathcal{O}_k \otimes_{A_{\mathfrak{p}}} \mathcal{J} 
\]
is an extension of the $\mathcal{E}_{\mathfrak{P}}$-module $\mathcal{J} = J \cap D_{\mathfrak{P}}$. And as above, any such $\mathcal{E}_{\mathfrak{P}}$-module is an $\mathcal{E}_{\mathfrak{P}}$-ideal, which we denote by $I_{\mathfrak{P}}$.\\

Collecting all the $\mathcal{I}_{\mathfrak{P}}$ for $\mathfrak{P} \neq \widetilde{\mathfrak{P}}$ to obtain $\mathcal{I}'_{\mathfrak{p}} = \oplus_{\mathfrak{P} \neq \widetilde{\mathfrak{P}}} \mathcal{I}_{\mathfrak{P}}$, we conclude that also every $(\mathcal{O}_k[f'_{\phi,\mathfrak{p}}] \otimes_{A_{\mathfrak{p}}} \mathcal{E}'_{\mathfrak{p}})$-submodule of $H^{\mathrm{\acute{e}t}}_{\mathfrak{p}}(\phi)$ is given by $\mathcal{I}'_{\mathfrak{p}}H^{\mathrm{\acute{e}t}}_{\mathfrak{p}}(\phi)$, as required.
 
This ends the proof of the theorem.
\end{proof}

\bibliographystyle{amsalpha}
\bibliography{Erratum.bib}

@misc{BKM,
      title={Abelian varieties over finite fields with commutative endomorphism algebra: theory and algorithms}, 
      author={Jonas Bergström and Valentijn Karemaker and Stefano Marseglia},
      year={2025},
      eprint={2409.08865},
      archivePrefix={arXiv},
      primaryClass={math.NT},
      url={https://arxiv.org/abs/2409.08865}, 
}

@article{BdV,
author = {de Vries, Sjoerd},
title = {Traces of {H}ecke operators on {D}rinfeld modular forms for $\mathrm{GL}_2(\mathbb{F}_q[T])$},
journal = {Mathematika},
volume = {72},
number = {2},
pages = {e70077},
doi = {https://doi.org/10.1112/mtk.70077},
year = {2026},
}

@article {GP2,
	AUTHOR = {Garai, Sumita and Papikian, Mihran},
	TITLE = {Computing endomorphism rings and {F}robenius matrices of
	{D}rinfeld modules},
	JOURNAL = {J. Number Theory},
	FJOURNAL = {Journal of Number Theory},
	VOLUME = {237},
	YEAR = {2022},
	PAGES = {pp.~145--164},
	ISSN = {0022-314X},
	DOI = {10.1016/j.jnt.2019.11.018},
	URL = {https://doi.org/10.1016/j.jnt.2019.11.018},
}

@book {Greither,
    AUTHOR = {Greither, Cornelius},
     TITLE = {Cyclic {G}alois extensions of commutative rings},
    SERIES = {Lecture Notes in Mathematics},
    VOLUME = {1534},
 PUBLISHER = {Springer-Verlag, Berlin},
      YEAR = {1992},
     PAGES = {x+145},
      ISBN = {3-540-56350-4},
       DOI = {10.1007/BFb0089165},
       URL = {https://doi.org/10.1007/BFb0089165},
}

@incollection {HayesCFT,
    AUTHOR = {Hayes, D.},
     TITLE = {Explicit class field theory in global function fields},
 BOOKTITLE = {Studies in algebra and number theory},
    SERIES = {Adv. in Math. Suppl. Stud.},
    VOLUME = {6},
     PAGES = {173--217},
 PUBLISHER = {Academic Press, New York-London},
      YEAR = {1979},
}

@article {KKP,
	AUTHOR = {Karemaker, Valentijn and Katen, Jeffrey and Papikian, Mihran},
	TITLE = {Isomorphism classes of {D}rinfeld modules over finite fields},
	JOURNAL = {J. Algebra},
	FJOURNAL = {Journal of Algebra},
	VOLUME = {644},
	YEAR = {2024},
	PAGES = {pp.~381--410},
}

@book {LaumonCDV,
    AUTHOR = {Laumon, G.},
     TITLE = {Cohomology of {D}rinfeld modular varieties. {P}art {I}},
    SERIES = {Cambridge Studies in Advanced Mathematics},
    VOLUME = {41},
 PUBLISHER = {Cambridge University Press, Cambridge},
      YEAR = {1996},
     PAGES = {xiv+344},
}

@book {DM,
	AUTHOR = {Papikian, Mihran},
	TITLE = {Drinfeld modules},
	SERIES = {Graduate Texts in Mathematics},
	VOLUME = {296},
	PUBLISHER = {Springer, Cham},
	YEAR = {2023},
	PAGES = {xxi+526},
	ISBN = {978-3-031-19706-2; 978-3-031-19707-9},
	DOI = {10.1007/978-3-031-19707-9},
	URL = {https://doi.org/10.1007/978-3-031-19707-9},
}

@article {Waterhouse,
	AUTHOR = {Waterhouse, William C.},
	TITLE = {Abelian varieties over finite fields},
	JOURNAL = {Ann. Sci. \'{E}cole Norm. Sup. (4)},
	FJOURNAL = {Annales Scientifiques de l'\'{E}cole Normale Sup\'{e}rieure. Quatri\`eme
	S\'{e}rie},
	VOLUME = {2},
	YEAR = {1969},
	PAGES = {pp.~521--560},
	ISSN = {0012-9593},
	URL = {http://www.numdam.org/item?id=ASENS_1969_4_2_4_521_0},
}

@article {JKYu,
	AUTHOR = {Yu, Jiu-Kang},
	TITLE = {Isogenies of {D}rinfeld modules over finite fields},
	JOURNAL = {J. Number Theory},
	FJOURNAL = {Journal of Number Theory},
	VOLUME = {54},
	YEAR = {1995},
	NUMBER = {1},
	PAGES = {pp.~161--171},
	ISSN = {0022-314X},
	DOI = {10.1006/jnth.1995.1108},
	URL = {https://doi.org/10.1006/jnth.1995.1108},
}

\end{document}